# Some Hodge theory from Lie algebras

Constantin Teleman[1]
teleman@math.stanford.edu

**Introduction**

The purpose of this note is two-fold. It is, first, a commentary on the (largely unknown) cohomology of Lie algebras of the form $\mathfrak{g} \otimes A$, where $\mathfrak{g}$ is a reductive Lie algebra and $A$ a finitely generated commutative $\mathbb{C}$-algebra; secondly, it briefly describes joint work with I. Grojnowski and S. Fishel [FGT], where this cohomology is computed in a special case, with noteworthy combinatorial applications. This relates to older work of Feigin's [F2].

The connection with Hodge theory is also two-fold. The theorem of Loday and Quillen [LQ] and Tsygan [Ts] describes the homology of $\mathfrak{gl}(A) := \mathfrak{gl} \otimes A$ as the free graded exterior algebra on the cyclic homology of $A$; this led Loday [L, Ch. X] to suggest that the homology of $\mathfrak{gl}_n(A)$ might be a free graded exterior algebra on a Hodge-style truncation of cyclic homology, defined from the *Adams* (or $\lambda$-) *decomposition*). Our computation confirms this in a special case. The bad news is that Loday's conjecture is *false* in general; the known examples seem to be 1-dimensional accidents, and our computation may well have exhausted the range of validity. However, a weaker conjecture of Feigin's [F1] for *smooth* algebras still stands.

In another direction, our Lie algebra computation determines the Hodge-to-de Rham (HdR) spectral sequence for the standard flag variety of a loop group $LG$ (reductive $G$). The failure of this spectral sequence to collapse (at $E_1$) was discovered by Grojnowki, who, with S. Fishel, computed its formal character, for $G = \mathrm{SL}_2$. (The present author is mainly responsible for the skepticism which led to delays in the project). For $G = \mathrm{SL}_\infty$, the $E_1$ term is easily computed from the theorem of Loday-Quilllen-Tsygan; but the general argument is more sophisticated. It relies on the collapse of the HdR spectral sequence of the moduli stack of holomorphic $G$-bundles over $\mathbb{P}^1$ [T1, §7], and on an older computation of Feigin's, reinterpreted in relation to the Hodge-to-de Rham sequence for the classifying stack $BG[z]$ of the group of $G$-valued polynomial maps.

**1. Refresher on cyclic homology**

We shall only discuss commutative or skew-commutative graded $\mathbb{C}$-algebras, following Burghelea-Vigué [BV] (cf. also [L, Ch. 5] for more details). The *cyclic homology* of $A$ has a bigrading

$$HC_n(A) = \bigoplus_{i=0}^{n} HC_n^{(i)}(A),$$

where $n$ is the homology degree and $i$ the Adams degree. Let first $A$ be smooth, all in degree zero, call $\Omega^k(A)$ the module of $k$th differentials over $A$, $d: \Omega^k \to \Omega^{k+1}$ de Rham's differential, and $H^k(A) := \mathrm{Ker}\, d/\mathrm{Im}\, d$ the $k$th de Rham cohomology of $\mathrm{Spec}(A)$. So $H^0(A)$ is the vector space of locally constant functions on $\mathrm{Spec}(A)$. $HC_n^{(i)}(A)$ is then described by the following table:

---

[1] Supported by an NSF Postdoctoral Fellowship

| $i \backslash n$ | 0 | 1 | 2 | 3 | 4 | 5 | 6 |
|---|---|---|---|---|---|---|---|
| 0 | $A$ | 0 | 0 | 0 | 0 | 0 | 0 |
| 1 | 0 | $\Omega^1/d\Omega^0$ | $H^0(A)$ | 0 | 0 | 0 | 0 |
| 2 | 0 | 0 | $\Omega^2/d\Omega^1$ | $H^1(A)$ | $H^0(A)$ | 0 | 0 |
| 3 | 0 | 0 | 0 | $\Omega^3/d\Omega^2$ | $H^2(A)$ | $H^1(A)$ | $H^0(A)$ |

Thus, $HC_\bullet^{(\bullet)}(\mathbb{C})$ contains one copy of $\mathbb{C}$ in each bi-degree $(n,i) = (2k,k)$, $k \geq 0$. In general, the $i$th row is the cohomology of the $i$-truncated de Rham complex $\Omega^{\leq i}$, flipped about degree $i$.

Take now a free skew-commutative algebra $A = \mathbb{C}[x_1, \ldots, x_N]$, with generators $x_i$ in degrees $d_i \geq 0$. (Skew-commutativity means $x_i x_j = (-1)^{d_i + d_j} x_j x_i$). The grading is homological: $A$ will soon carry a differential of degree $(-1)$. Its de Rham complex is the free skew-commutative algebra
$\Omega_\bullet^\bullet(A) := A[dx_1, \ldots, dx_n]$, $dx_i$ being a free generator of degree $d_i - 1$. (So $(dx)^2 \neq 0$ if $x$ is odd). De Rham's operator $d: x_i \mapsto dx_i$ extends to a graded derivation of (cohomological) degree 1. The superscript grading in $\Omega_\bullet^\bullet(A)$ is the degree of the form (the number of $dx_i$ factors), whereas the subscript grading (the algebra degree) counts both $x_i$ and $dx_i$ with weight $d_i$. The difference (form degree – algebra degree) is the *total degree*. This is the usual grading, if all $d_i$ are nil. This time, $HC_\bullet^{(i)}(A)$ is the cohomology of the truncated de Rham complex $\Omega_\bullet^{\leq i}(A)$, truncated by the form degree, but graded by the *total* degree, again flipped about $i$ (so $HC_n^{(i)}$ corresponds to total degree $2i - n$).

Let now $A$ carry a graded derivation $\partial$ of degree $(-1)$, with $\partial^2 = 0$: $(A, \partial)$ is then a commutative DGA. $\partial$ extends to $\Omega_\bullet^\bullet(A)$, in such a way that $(d + \partial)^2 = d\partial + \partial d = 0$. With the differential $d + \partial$, of degree $(+1)$ for the total grading, $\Omega_\bullet^\bullet(A)$ is the total de Rham complex of $(A, \partial)$. Define $HC_n^{(i)}(A, \partial)$ to be the cohomology, in total degree $(2i - n)$, of its $i$-truncation $\Omega_\bullet^{\leq i}(A)$. One checks that a quasiisomorphism $f: (A, \partial) \to (A', \partial')$ (a morphism of DGA's inducing $\partial$-homology isomorphisms) gives an isomorphism in cyclic homology. Any commutative $A$ can be resolved by a free DGA $A \xleftarrow{\varepsilon} A_0 \leftarrow A_1 \leftarrow A_3 \leftarrow \ldots$, so this allows one to define the cyclic homology of commutative (or graded-commutative) algebras.

*Remarks.* (i) Over a smooth proper scheme $X$, the sheafified cyclic homology is a Hodge truncation of de Rham cohomology: $HC_n^{(i)}(X) = H^{2i-n}(X)/F^{i+1}H^{2i-n}(X)$ [W].

(ii) $\Omega_\bullet^{\leq i}$ is a quotient complex of $\Omega_\bullet^\bullet$, so a natural *periodicity map* $S: HC_n^{(i)}(A) \to HC_{n-2}^{(i-1)}(A)$ is defined. Feigin and Tsygan [FT2] show that $\lim_k HC_{2k-n}^{(k)}(A)$ is the $n$th crystalline cohomology of $\text{Spec}(A)$.

## 2. Homology of $\mathfrak{gl}(A)$

Call $\mathfrak{gl}(A)$ the complex Lie algebra of finite rank, infinite matrices with entries in $A$, with the commutator bracket. (If $A$ is graded, using graded commutators defines a graded Lie algebra; in general, for any complex Lie algebra $\mathfrak{g}$ and commutative DGA $A$, there is a natural structure of differential graded Lie algebra on $\mathfrak{g}(A) := \mathfrak{g} \otimes_\mathbb{C} A$.) Lie algebra homology is taken with $\mathbb{C}$ coefficients, unless others are shown.

**(2.1) Theorem.** (Loday-Quillen [LQ], Tsygan [Ts]) $H_\bullet(\mathfrak{gl}(A)) = \mathbb{C}[HC_\bullet(A)[-1]]$.

The right-hand side is the free graded-commutative algebra generated by the graded vector space with $k$th piece $HC_{k-1}(A)$. The isomorphism preserves the graded Hopf algebra structures, using the natural comultiplication on Lie algebra homology, and the multiplication induced by the "Whitney sum" of matrices. The theorem only requires associativity of $A$; but we are also interested in other $\mathfrak{g}$'s, for which things seem hopeless, unless $A$ is commutative, or at least graded-commutative. The map in (2.1) is easily described when $A$ is smooth. For any $\mathfrak{g}$ and a homogeneous invariant polynomial $P$ of degree $(i+1)$ on $\mathfrak{g}$, consider for $i \le n \le 2i$ the linear maps $\varphi_{n,P}: \wedge^{n+1}\mathfrak{g}(A) \to \Omega^{2i-n}(A)$

(2.2) $\qquad \varphi_{n,P}(\xi_0 \wedge \ldots \wedge \xi_n) := \mathrm{Alt}\, P(\xi_0, [\xi_1, \xi_2], \ldots, [\xi_{2n-2i-1}, \xi_{2n-2i}], d\xi_{2n-2i+1}, \ldots, d\xi_n),$

"Alt" being the complete antisymmetrization. Direct computation verifies the following ([F1], but given, for some reason, without the formula). $\partial$ is the Lie algebra differential in $\wedge^\bullet \mathfrak{g}(A)$:

**(2.3) Proposition.** *(a)* $\varphi_{n,P} \circ \partial = \mathrm{const} \cdot d \varphi_{n+1,P}$ *if* $i \le n < 2i$; $\varphi_{2i,P} \circ \partial = 0$.
*(b)* $\varphi_{n,P}$ *takes Lie algebra boundaries to exact forms, and, if* $i < n \le 2i$, *it takes cycles to closed forms.*
*(c) The* $\varphi_{n,P}$ *assemble to a linear map* $\varphi_P : H_{\bullet+1}(\mathfrak{g}(A)) \to HC_\bullet^{(i)}$.
*(d) When* $\mathfrak{g} = \mathfrak{gl}$ *and P ranges over the* basic *invariant polynomials* $M \mapsto \mathrm{Tr}(M^{i+1})$, $i \ge 0$, *the* $\varphi_P$ *map the space of primitive elements in* $H_\bullet(\mathfrak{gl}(A))$ *isomorphically onto* $HC_\bullet(A)$.

Readers who prefer algebras to coalgebras can think of the dual map to (c) geometrically, as follows: the integral of $\varphi_P$ over a compact $(2i-n)$-cycle $C$ on $\mathrm{Spec}(A)$ is a Lie algebra $(n+1)$-cocycle, whose class depends only on the homology class of $C$, if $n > i$. Calling $\delta$ the Lie algebra cohomology differential,

(2.4) $\qquad \delta \int_C \varphi_n^P = \int_C \varphi_n^P \circ \partial = \mathrm{const} \int_C d \varphi_{n+1}^P = \mathrm{const} \int_{\partial C} \varphi_{n+1}^P = 0.$

*(2.5) Remark.* The proposition extends to smooth DGA's, after suitable adjustment in the grading; and this describes the isomorphism (2.1) for all commutative algebras.

### 3. Conjectures and some results for $\mathfrak{g}(A)$

The L-Q-T theorem (2.3.d) begs to be extended to other Lie algebras. For reductive $\mathfrak{g}$, there is an obvious candidate, conjectured, in slightly different forms, by Feigin [F1] and Loday [L]. Recall that $(\mathrm{Sym}^\bullet \mathfrak{g}^t)^G$ is a polynomial algebra with generators of degrees $m_1+1, \ldots, m_l+1$, where $l$ is the rank and $\{m_i\}$ are the exponents of $\mathfrak{g}$; (2.3.c) defines a homomorphism of graded-commutative algebras

(3.1) $\qquad \Phi: \mathbb{C}\left[\bigoplus_{i=1}^l HC_\bullet^{(m_i)}(A)^t[-1]\right] \to H^\bullet(\mathfrak{g}(A)).$

When $\mathfrak{g} = \mathfrak{gl}_n$, the image of $\Phi$ is that of $H^\bullet(\mathfrak{gl}(A)) \to H^\bullet(\mathfrak{gl}_n(A))$ under the natural restriction, and the computation of the stable range in [LQ] shows that $\Phi$ is an isomorphism in degrees no higher than $n$. Also, when $A = \mathbb{C}$ and any reductive $\mathfrak{g}$, (3.1) is the well-known description of $H^\bullet(\mathfrak{g})$ as a free exterior algebra. Lie algebra cocycles give left-invariant forms on the group, identifying $H^\bullet(\mathfrak{g})$ with $H^\bullet(G; \mathbb{C})$; and (3.1) corresponds to the topological fact that $G$ is, in rational homotopy, a product of odd spheres of dimensions $2m_i + 1$.

Remarkably, turns out to be an *isomorphism* in a number of non-trivial cases. Here are the examples I know, with one to be added at the end of this section (cf. also [L, 10.3.10]):

(i)   $A = \mathbb{C}[x]$: "classical", from the BGG resolution; cf. [G], [F1].
(ii)  $A = C^\infty(S^1)$: [PS], [F1] (with continuous cochains); [K2], for $A = \mathbb{C}[x, x^{-1}]$.
(iii) $A = \mathbb{C}[\ ]$, any smooth affine curve: [T2] (but probably "known to experts")
(iv)  $A = \mathbb{C}[x,\xi]$, the graded algebra in one even and one odd variable: [F2].

The non-zero cyclic homologies are (for $i > 0$ since $HC^{(0)} = A$, always)

(3.2) $$HC_n^{(i)}(\mathbb{C}[x,\xi]) = \begin{cases} \mathbb{C} \oplus \mathbb{C}[x]dx\ \xi(d\xi)^{i-1}, & n = 2i \\ \mathbb{C}[x]\ \xi(d\xi)^i, & n = 2i+1 \end{cases}$$

*(3.3) Remark.* In more conventional terms, $H_q(\mathfrak{g}[x,\xi])$ equals $\bigoplus_p H_{q-2p}(\mathfrak{g}[x]; \mathrm{Sym}^p\{\xi \mathfrak{g}[x]\})$.

(v)  $A = \mathbb{C}[x]/(x^n)$: conjectured (and proved for $\mathfrak{sl}_n$) by Hanlon [H]; ingeniously proved by Feigin [F2] in general, using the previous computation (iv). To describe the cyclic homology, resolve $A$ by $(\mathbb{C}[x,\xi],\ )$, with $\xi = x^n$. We get from (3.2), keeping track of the $x$-grading,
$HC_n^{(i)}(\mathbb{C}[x]/(x^n)) = \mathbb{C} \oplus x^{ni+1}\mathbb{C}[x]/x^{n(i+1)}\mathbb{C}[x]$ if $n = 2i$, zero otherwise.

(vi) $A = \mathbb{C}[\ ][\xi]$, a smooth curve and $\xi$ an odd variable [FGT]. Generalizing (3.2),

(3.4) $$HC_n^{(i)}(\mathbb{C}[x,\xi]) = \begin{cases} H^1(\ ; \mathbb{C}), & n = 2i-1 \\ \mathbb{C} \oplus \Omega^1[\ ]\xi(d\xi)^{i-1}, & n = 2i \\ \Omega^0[\ ]\xi(d\xi)^i, & n = 2i+1 \end{cases}$$

For reference, here are the Lie algebra homologies for *semi-simple* $\mathfrak{g}$. For reductive $\mathfrak{g}$, the center splits off, contributing factors of $\Lambda(A)$ to the homology of $\mathfrak{g}(A)$.

$H_*(\mathfrak{g}[x]) \cong H_*(\mathfrak{g})$;
$H_*(C^\infty(S^1;\mathfrak{g})) \cong H_*(C^\infty(S^1;G);\mathbb{C})$;
$H_*(\mathfrak{g}[\ ]) = H_*(C^\infty(\ ;G);\mathbb{C})$;
$H_*(\mathfrak{g}[x,\xi])$ is freely generated, over $H_*(\mathfrak{g})$, by copies of $\mathbb{C}[x]dx$ and $\mathbb{C}[x]$, in degrees $2m_i + 1$ and $2m_i + 2$, respectively;
$H_*(\mathfrak{g}[x]/(x^n)) \cong H_*(\mathfrak{g})^{\otimes n}$ (as a coalgebra, but non-canonically).

*(3.5) Remark.* Feigin's proof [F2] of (iv) contains an error for which I know no cure (the key Lemma on p.93 is wrong in absolute homology, where $Q = \partial/\partial\xi$ induces a non-zero differential $d_2: E^2_{2,2} \to E^2_{0,3}$; and the argument does not seem to distinguish between absolute and relative homology). A different proof will be given in [FGT].

Loday [L, Ch.10] made the beautiful suggestion that might always be an isomorphism, at least when $\mathfrak{g} = \mathfrak{gl}_n$. This would be a far-reaching strengthening of the rational homotopy splitting of $GL_\infty$ into $K(\mathbb{Q};2i+1)$'s. Unfortunately, this fails in the simplest 2-dimensional example, $\mathfrak{sl}_2[x,y]$ (and then, of course, also for $\mathfrak{gl}_2 = \mathfrak{gl}_1 \oplus \mathfrak{sl}_2$).

**(3.6) Proposition.** *The square of the class* $\xi \otimes \eta \mapsto \mathrm{Tr}(\xi_x \eta_y - \xi_y \eta_x)$ *is null in* $H^4(\mathfrak{sl}_2[x,y])$.

Indeed, $H(\mathfrak{sl}_2[x,y]) = H(\mathfrak{sl}_2) \otimes H(\mathfrak{sl}_2[x,y], \mathfrak{sl}_2)$, and the second factor is computed by the $\mathfrak{sl}_2$-invariant part of the chain complex $(\mathfrak{sl}_2[x,y]/\mathfrak{sl}_2)$. Counting dimensions of the invariant spaces of $(x,y)$-type $(2,2)$ gives, starting in degree zero, $0,0,3,3,1,0,0,\ldots$. As there is one relative class in $H_2$ (stable range) and none in $H_3$ [C], the 4-cycle must be a boundary. So there is no class of the correct type in $H^4$. Note that the 2-cocycle in (3.6) picks out the coefficient of $dx \wedge dy$ in $H_2 = \Omega^1/d\Omega^0$, so it is not the zero class.

Feigin [F1] proposed a more prudent conjecture: for smooth $A$, $\chi$ should be *onto*. For $\mathfrak{gl}_n$, this says that all classes come from $\mathfrak{gl}$. Without smoothness, this can fail, as seen in the case of two crossing lines, $A = \mathbb{C}[x,y]/(xy)$: $HC_n^{(1)}(A)$ is $\mathbb{C}$ in degrees $n = 1,2$, null otherwise; and, while $H(\mathfrak{sl}_2(A))$ indeed looks additively like a free graded algebra, with generators in degrees 2 and 3, the cup-product $H^2 \otimes H^2 \to H^4$ is nil [T2]. Thus, the generator of $H^4(\mathfrak{sl}_2(A))$ cannot come from $\mathfrak{sl}$. (This can be explained by viewing $A$ as a degeneration of $\mathbb{C}[x,x^{-1}]$: additively, $H(\mathfrak{sl}_2(A))$ remains unchanged, but the ring structure degenerates. Note that $HC_1^{(1)}(A)$ sees the vanishing cycle). To my knowledge, Feigin's conjecture still stands, and even the case $A = \mathbb{C}[x,y]$ is open. I am not aware of a standing conjecture for the relations in the cohomology ring.

The last example where $\chi$ is an isomorphism concerns a first-order extension of $\mathbb{C}[x]$:

(vii)  $A = \mathbb{C}[x] \oplus M$, where $M$ is the $\mathbb{C}[x]$-module $\mathbb{C}[x,x^{-1}]/\mathbb{C}[x]$.

Its Lie algebra homology decomposes $H_p(\mathfrak{g}(A)) = \bigoplus_q H_{p-q}(\mathfrak{g}[x]; \wedge^q M)$; cf. (3.3), and has a noteworthy Hodge-theoretic interpretation, as we shall see in the next section. Note that $\chi$ is no longer iso when $M = \mathbb{C}[x]$, e.g. for $\mathfrak{sl}_2(A) = \mathfrak{sl}_2[x,y]/(y^2)$, for the same reason as (3.6) (there is a class in $H^2$ squaring to zero). We are skirting the bounds of validity Loday's conjecture.

The non-zero cyclic homology groups $HC_n^{(i)}$ in question are (again, $i > 0$) $HC_{2i-1}^{(i)} = M \otimes dx$ and $HC_{2i}^{(i)} = \mathbb{C} \oplus M$. We shall show in [FGT] that the relative homology $H_q(\mathfrak{g}[x], \mathfrak{g}; \wedge^p\{\mathfrak{g}[x,x^{-1}]/\mathfrak{g}[x]\})$ is the free bigraded coalgebra on copies of $M \otimes dx$, in bi-degrees $(p,q) = (m_i, m_i)$ and $M$, in bi-degrees $(m_i + 1, m_i)$ cf. (5.3) and (5.5) below.

## 4. Loop groups

The last Lie algebra homology arose naturally in [FGT] as the (dual of the) $E_1$ term in Hodge-to-de Rham spectral sequence for the flag variety of the loop group of $G$. This Hodge-theoretic link seems unrelated to cyclic homology; but the latter offers a comforting explanation for our answers.

For a simple Lie group $G$, let $LG := G[z,z^{-1}]$ be the (algebraic ind-)group of $G$-valued Laurent polynomials. (Bounding the order of the poles at $0$ and $\infty$ realizes $LG$ as a union of closed, finite-dimensional subvarieties). The homogeneous (ind-)variety $X := G[z,z^{-1}]/G[z]$ carries an ample line bundle $\mathcal{L}$, generator of $\mathrm{Pic}(X) \cong \mathbb{Z}$, a suitable power of which gives a closed Kodaira embedding of $X$ in $\mathbb{P}(\mathcal{H})$ [K], [M]; $\mathcal{H}$ is, of course, dual to the space of sections $(X;\mathcal{L}^h)$. is filtered downwards by $G[z]$-invariant subspaces of finite codimension, while $\mathcal{H}$ is increasingly filtered by finite-dimensional spaces. Both filtrations are associated to an *energy* grading, positive on , negative on $\mathcal{H}$, which comes from the $\mathbb{C}^\times$-action on $(X;\mathcal{L}^h)$ defined by $z$-rotation. The filtration on $\mathcal{H}$ gives an ind-variety structure on $\mathbb{P}(\mathcal{H})$ compatible with the ind-structure on $X$. Thus, $X$ is a union of projective varieties. It is also formally smooth and reduced [LS], and can be covered by Zariski-open subsets isomorphic to $G[z^{-1}]/G$.

Another version of these objects are the *formal loop group* $\hat{L}G := G((z^{-1}))$ and its *thick flag variety* $\hat{X} := G((z^{-1}))/G[z]$. This is the closure of $X$ in $\mathbb{P}(\hat{\mathcal{H}})$, where $\hat{\mathcal{H}}$ is the direct product of the energy eigenspaces in $\mathcal{H}$, and is a scheme of infinite type, covered by affine spaces isomorphic to $G[[z^{-1}]]/G$. The thick flag variety has the virtue of being smooth, while the thin one is a union of compact varieties. (This "ind-compactness" was used by Kumar and Mathieu in their proofs of the Borel-Weil-Bott theorem). For some purposes, the two flag varieties are interchangeable; often, cohomologies of equivariant vector bundles (under loop rotation) on the two are related in an obvious way, the one over $X$ being a direct product, the one over $\hat{X}$ a direct sum, of the same $\mathbb{C}^\times$-eigenspaces.

The geometric study of highest-weight representations of $LG$ relies on the philosophy that $X$ and $\hat{X}$ behave, in most relevant ways, like a compact Kähler manifold, being analogous to the flag varieties of compact Lie groups. The failure of their HdR spectral sequence to collapse at $E_1$ was thus unexpected. Before describing the sequence, let us note one geometric consequence.

**(4.1) Proposition.** *$X$ is not a complex manifold, and $G[z]$ is not a complex Lie group. (As analytic ind-varieties, they are not locally isomorphic to vector spaces).*

Roughly speaking, if $X$ were a manifold, its de Rham complex with the naïve filtration would be a pure Hodge complex, and "ind-compactness" would force the collapse at $E_1$ of HdR. (Making this more precise requires a bit of technology, and will wait for [ST]). So every point of these homogeneous spaces is equally singular. Recall again that $X$ is reduced and formally smooth; for schemes (of finite type), this implies smoothness, so we are seeing here a bizarre infinite-dimensional phenomenon. ($\hat{X}$, on the other hand, is smooth, being covered by affine spaces; but it is not covered by compact varieties).

## 5. Hodge-to de Rham spectral sequence for $X$

For simplicity, $G = \mathrm{SL}_2$; similar statements hold for any semi-simple $G$. Since $X$ is homotopy equivalent to the smooth based loop group $G$ [PS, Ch. 9], $H^*(X;\mathbb{C})$ is the symmetric algebra generated by the Chern class $c_1(\mathcal{L})$. This has Hodge type $(1,1)$, and this determines the $(p,q)$ decomposition on $H^*(X)$. (In fact, the Hodge structure on $H^*(X)$ is pure, of Tate type, for any $G$: the algebra generators transgress from those of $H^*(BG)$, whose Hodge structure was described by Deligne). Incidentally, $c_1(\mathcal{L})$ has a representative that is invariant under the translation action of $LG^{\mathrm{cpt}}$,

(5.1) $$\omega(\xi,\bar{\eta}) = \frac{1}{2} \circ \langle \xi | d\bar{\eta} \rangle$$

for two (holomorphic) tangent vectors $\xi, \eta \in \mathfrak{g}[z^{-1}]/\mathfrak{g} = T^{1,0}X$ at the base-point of $X$.

The Hodge-to-de Rham spectral sequence for $X$ is

(5.2) $$E_1^{p,q} = H^q(X; \Omega^p) \Rightarrow H^{p+q}(X; \mathbb{C}).$$

As a vector bundle, $\Omega^1 = G[z, z^{-1}] \times_{G[z]} \mathfrak{g}[[z]]dz$, and general theory gives canonical identifications

(5.3) $$H^q(X; \Omega^p) = H^q_{G[z]}\left( \Lambda^p \{\mathfrak{g}[[z]]\, dz\} \right) = H^q(\mathfrak{g}[z], \mathfrak{g}; \Lambda^p \{\mathfrak{g}[[z]]\, dz\}),$$

of the desired spaces as either algebraic group cohomologies or Lie algebra cohomologies.

*(5.4) Remark.* One step in (5.3) that is worth mentioning is that $H^q(X; \Omega^p)$ is trivial, as a representation of $LG$. For the flag variety of a compact Lie group, this space is a summand in $H^{p+q}(X; \mathbb{C})$, on which the (connected) group must act trivially. This argument breaks down in our case, but we can use the general fact that the only genuine (non-projective) highest-weight representation of $LG$ is the trivial one.

Were (5.2) to collapse at $E_1$, we would have $H^q(X; \Omega^p) = \mathbb{C}$, if $q = p$, and zero otherwise. Let us see why this cannot be true. Consider $p = q = 1$. One readily finds a 1-cocycle of $G[z]$ with values in $\mathfrak{g}[[z]]$, simply $\gamma \mapsto d\gamma \cdot \gamma^{-1}$ (or rather, its formal germ at the origin). This can be traced back to the Kähler class (5.1), as the derivative suggests. However, the cocycle can be multiplied by any scalar-valued formal power series, and yields more non-trivial group cohomology classes. These are topologically unaccounted for, and must transgress under some differentials in (5.2).

It turns out there are no other 1-classes: $E_1^{1,1} = \mathbb{C}[[z]]$. One can easily compute a bit more. For instance, $E_1^{0,q} = 0$ if $q > 0$, from our knowledge of $H^\bullet(\mathfrak{g}[x]) \cong H^\bullet(\mathfrak{g}[x], \mathfrak{g}) \otimes H^\bullet(\mathfrak{g})$; also, $E_1^{p,0} = 0$ if $p > 0$ (no invariants). It is also known [G] that $E_1^{1,q} = 0$ for $q > 1$. In the simplest of all worlds, the first unknown group $E_1^{2,1}$ would be the $d_1$-image of the non-surviving part $\mathbb{C}[[z]]/\mathbb{C}$ of $E_1^{1,1}$, which we identify with $\mathbb{C}[[z]]\, dz$, by the de Rham differential. The optimistic guess is correct.

**(5.5) Theorem.** [FGT] *The $E_1$ term in (5.2) is the bigraded skew-commutative algebra freely generated by $\mathbb{C}[[z]]$, in bi-degree $(1,1)$, and $\mathbb{C}[[z]]dz$, in bi-degree $(1,2)$: $\mathrm{Sym}^\bullet \mathbb{C}[[z]] \otimes \Lambda^\bullet \mathbb{C}[[z]]dz$; and $d_1$ is the graded algebra derivation generated by de Rham's operator $d : \mathbb{C}[[z]] \to \mathbb{C}[[z]]dz$. The sequence collapses at $E_2$.*

The proof would take us to far afield and will wait for [FGT]. Instead, I shall briefly describe the calculation [Gr] that led to the discovery of (5.5), at a time when the three of us were blissfully ignorant of the contents of the present §2–3.

While there is no "textbook" method for computing $H^n(X; \Omega^m)$, a formula for the holomorphic Euler characteristic $\chi(V) = \sum_n (-1)^n H^n(X; V)$ of homogeneous vector bundles $V$ exists; it combines the Kac character formula for integrable representations of $LG$ with (a weak form of) the Borel-Weil-Bott formula of Kumar and Mathieu. The formula gives the character of the natural $G \times \mathbb{C}^\times$-action on the vir-tual space $\chi(V)$, where $G \subset LG$ acts by translation and $\mathbb{C}^\times$ by $z$-rotation; we assume that the rotation action lifts holomorphically to $V$. Subject to some restrictions on $V$, this virtual character is a formal power series $\mathrm{ch}(X;V)(u,q) := \mathrm{Tr}_{\chi(V)}(u\ q) \in \mathbb{Z}[u,u^{-1}]((q))$, where $u$ is in a fixed maximal torus $\mathbb{C}^\times \subset SL_2$ and $q \in \mathbb{C}^\times$ (the rotations). (This reuse of $q$, compared to (5.2), is unfortunate, but both are traditional).

Note that $G \times \mathbb{C}^\times$ acts on the fibre $V_1$ of $V$ over $1 \in X$; call its character $\mathrm{ch}V_1(u,q)$. Two more notational elements are needed: the action of the *affine Weyl group* $\mathbb{Z} \tilde{\times} \mathbb{Z}/2$ on the ring $\mathbb{Z}[u,u^{-1}]((q))$, under which $n \in \mathbb{Z}$ sends $u$ to $uq^n$, while $(-1) \in \mathbb{Z}/2$ sends $u$ to $u^{-1}$ (both preserve $q$), and the convention, standard in the theory of hypergeometric series,

$$(5.6) \qquad (a)_n = \begin{cases} (1-a)(1-qa)\cdots(1-q^{n-1}a), & 0 \le n \\ (1-q^na)^{-1}\cdots(1-q^{-2}a)^{-1}(1-q^{-1}a)^{-1}, & n < 0 \end{cases}$$

Finally, the Kac formula for the holomorphic Euler characteristic of $V$ is

$$(5.7) \qquad \mathrm{ch}(X;V)(u,q) = \sum_{w \in \mathbb{Z}\tilde{\times}\mathbb{Z}/2} w \cdot \frac{\mathrm{ch}V_1(u,q)}{(u^{-2})_\infty (q)_\infty (qu^2)_\infty}$$

The $\mathrm{ch}(X;\wedge^m)(u,q)$ are best captured by their Poicaré series $\sum_m (-t)^m \mathrm{ch}(X;\wedge^m)(u,q)$. The fiber $V_1$ over $1 \in X$ is $(\mathfrak{g}[[z]] \ dz)$, cf. (5.3); its character is $\prod_{n>0}(1-tq^n u^{-2})(1-tq^n)(1-tq^n u^2)$, if we keep track of the powers of $t$. We obtain from (5.7)

$$(5.8) \qquad \sum_{m \ge 0} (-t)^m \mathrm{ch}(X;\wedge^m)(u,q) = \sum_{w \in \mathbb{Z}\tilde{\times}\mathbb{Z}/2} w \cdot \frac{(tqu^{-2})_\infty (tq)_\infty (tqu^2)_\infty}{(u^{-2})_\infty (q)_\infty (qu^2)_\infty},$$

which, after unraveling the $w$-action, is converted by elementary manipulations into

$$(5.9) \qquad \sum_{m \ge 0} (-t)^m \mathrm{ch}(X;\wedge^m)(u,q) = \frac{(tqu^{-2})_\infty (tq)_\infty (tqu^2)_\infty}{(qu^{-2})_\infty (q)_\infty (qu^2)_\infty} \sum_{n \in \mathbb{Z}} \frac{1-q^{2n}u^2}{1-u^2} \frac{(t^{-1}u^2)_{2n}}{(tqu^2)_{2n}} t^{2n}.$$

Remarkably, a summation formula for (5.9) exists. In slightly more general form, it is known as *Ramanujan's sum for the* $_1\psi_1$ *series* and reads

$$(5.10) \qquad \sum_{n=-\infty}^{\infty} \frac{(a)_n}{(b)_n} t^n = \frac{(q)_\infty (b/a)_\infty (at)_\infty (q/at)_\infty}{(b)_\infty (q/a)_\infty (t)_\infty (b/at)_\infty}.$$

We leave to the reader the joy of performing the variable substitutions in (5.10) that reduce (5.9) to

$$(5.11) \qquad \sum_{m,n} (-1)^n (-t)^m \mathrm{ch} H^n(X;\wedge^m)(u,q) = \frac{1}{1-t} \frac{(qt^2)_\infty}{(qt)_\infty},$$

the final answer. The numerator on the right is the contribution of $\mathbb{C}[[z]]dz$ from (5.5) (generated in $m=2, n=1$), the denominator that of $\mathrm{Sym}\,\mathbb{C}[[z]]$ ($m=n=1$). The topological part is $1/(1-t)$.

For other Lie groups, the identities arising by equating the analogue of (5.9) with our answer

(5.5) seem new. At any rate, combinatorics alone would not compute the spectral sequence, since the $H^\bullet(X; {}^m)$ are not pure-dimensional. Instead, the proof of (5.5) goes roughly as follows: $X$ is the fiber of a natural morphism from the stack $\mathfrak{M}$ of $G$-bundles over $\mathbb{P}^1$ to the classifying stack $BG[z^{-1}]$. This fibration follows from the *uniformization theorem* [LS], which equates $\mathfrak{M}$ with the double coset stack $G[z^{-1}] \setminus G[z,z^{-1}]/G[z]$. In [T1], it is shown that the HdR spectral sequence for $\mathfrak{M}$ collapses at $E_1$; this determines $H^n(\mathfrak{M}; {}^m)$. On the other hand, the HdR sequence for $BG[z^{-1}]$ is closely related to Feigin's calculation (iv). From the two and from the Leray spectral sequence, $H^\bullet(X; {}^m)$ can be extracted.